\documentclass{amsart}
\usepackage{graphicx}
\usepackage{float}
\usepackage[arc,all]{xy}
\usepackage{amsthm}
\newtheorem{thm}{Theorem}[section]
\newtheorem{cor}[thm]{Corollary}
\newtheorem{lem}[thm]{Lemma}

\theoremstyle{definition}

\newtheorem{exa}[thm]{Example}
\theoremstyle{remark}

\numberwithin{equation}{section}
\newtheorem*{theorem}{\textbf{Main Theorem}}

\begin{document}

\title{ On the dimension of Schur multiplier of  Lie algebras of maximal class }
\author[A. Shamsaki]{Afsaneh Shamsaki}
\address{School of Mathematics and Computer Science\\
Damghan University, Damghan, Iran}
\email{Shamsaki.Afsaneh@yahoo.com}
\author[P. Niroomand]{Peyman Niroomand}
\email{niroomand@du.ac.ir, p$\_$niroomand@yahoo.com}
\address{School of Mathematics and Computer Science\\
Damghan University, Damghan, Iran}


\begin{abstract}
The paper is devoted to obtain an upper bound for the Schur multiplier of nilpotent Lie algebras of maximal class. It improves the later upper bounds on the Schur multiplier of such Lie algebras.
\end{abstract}
\maketitle
\section{Introduction and Preliminaries}
There are wide literature to show that the results on the Schur multiplier of finite $ p$-group have analogues on the Schur multiplier $\mathcal{M}(L)$ of a nilpotent Lie algebra $ L $. \\
For an $ n $-dimensional nilpotent Lie algebra,  Moneyhun in \cite{2} showed $ \dim \mathcal{M}(L)\leq \frac{1}{2}n(n-1) $. Later this bound improved by the authors in \cite{3}. More precisely, it is proved for a non-abelian nilpotent Lie algebra of dimension $ n $ with the derived subalgebra of dimension $ m $, $ \dim \mathcal{M}(L)\leq \frac{1}{2}(n+m-2)(n-m-1)+1 $. When $ L $ is nilpotent of maximal class, then $ \dim L^{2}=n-2 $ and so we have $ \dim \mathcal{M}(L)\leq n-1 $. Recently this bound improved for the of nilpotent Lie algebra of maximal class in  $[4$, Theorem $3.1]$. It is shown that for maximal class for a nilpotent Lie algebra of maximal class of dimension $ n $, $ \dim \mathcal{M}(L)\leq n-2 $. \\
In this paper for a field of characteristic not equal to  2, we are going to improve this bound as follows.
\begin{theorem}
{\em Let $ L $ be an $ n $-dimensional nilpotent Lie algebra of maximal class. Then
$\dim \mathcal{M}(L)\leq n/2$ when $ n $ is even and $\dim \mathcal{M}(L)\leq \lceil(n+1)/2\rceil$ otherwise.}
\end{theorem}
 We also give some examples of  nilpotent Lie algebras of maximal class such that the Schur multiplier of them obtained the upper bound of Main Theorem. (See Examples $ 2.1 $ and $ 2.2 $.)\\
Let $ F/R $ be a free presentation of the $ c $-step nilpotent Lie algebra $ L $ and let for all $ 2\leq i \leq c $,
$\lambda_{i}: L/\gamma_{2}(L) \otimes \gamma_{i}(L)/\gamma_{i+1}(L) \longrightarrow [F, \gamma_{i}(F)+R]/[\gamma_{i+1}(F)+R, F]$  are defined  in  \cite{5}. Then
\begin{lem}$($See $[5$, Lemma $3.1])$
Let $ L $ be a Lie algebra. Then for all $ x_{1},x_{2},...,x_{i+1}\in L $, we have
\begin{align*}
&[[x_{1},x_{2},...,x_{i}]_{l},x_{i+1}]+[[x_{i+1},[x_{1},x_{2},...,x_{i-1}]_{l}],x_{i}]\cr
&+[[[x_{i},x_{i+1}]_{r},[x_{1},...,x_{i-2}]_{l}],x_{i+1}]\cr
&+[[[x_{i-1},x_{i},x_{i+1}]_{r},[x_{1},x_{2},...,x_{i-3}]_{l}],x_{i-2}]+...+[[x_{2},...,x_{i+1}]_{r},x_{1}]=0
\end{align*}
for $ i\geq 3 $, where
$$  [x_{1},x_{2},...,x_{i}]_{r}=[x_{1},[...[x_{i-2},[x_{i-1},x_{i}]]...]  ~  and   ~   [x_{1},x_{2},...,x_{i}]_{l}=[...[[x_{1},x_{2}],x_{3}],...,x_{i}].      $$

\end{lem}
 \begin{cor} $($See $[5$, Corollary $3.2])$
 Let $ \Psi_{i} : L\times ... \times L \longrightarrow \gamma_{i}(L)/\gamma_{i+1}(L)\otimes L/\gamma_{2}(L)$, given by  \begin{align*}
\Psi_{i}(x_{1}, x_{2},...,x_{i+1}) =&\overline{[x_{1}, x_{2},...,x_{i}]}\otimes \overline{x_{i+1}}+\overline{[x_{i+1}[x_{1}, x_{2},...,x_{i-1}]_{l}]}\otimes \overline{x_{i}}\cr
& \overline{[[x_{i},x_{i-1}]_{r},[x_{1}, x_{2},...,x_{i-2}]]}\otimes \overline{x_{i-1}}\cr
& \overline{[[x_{i-1},x_{i},x_{i+1}]_{r},[x_{1}, x_{2},...,x_{i-3}]_{l}]}\otimes \overline{x_{i-2}}\cr
& +...+ \overline{[x_{2},...,x_{i+1}]}\otimes \overline{x_{1}},
\end{align*}
for all $ 2\leq i \leq c $. Then $\mathrm{Im}\Psi_{i} \subseteq  \ker \lambda_{i}$.
\end{cor}
To prove the Main Theorem we need the following theorem.
\begin{thm}$($See $[3$, Corollary $2.3])$
Let $L$ be a finite dimensional Lie algebra and $K$ be a central ideal of $L$. Then
\begin{align*}
\dim\mathcal{M}(L)+\dim(L^{2}\cap K) \leq \dim\mathcal{M}(L/K)+\dim\mathcal{M}(K)+\dim (L/K)^{ab}\otimes K.
\end{align*}
\end{thm}
\section{proof of the main theorem}
 {\bf Proof of the Main Theorem.} Since $ L $ is an $ n $-dimensional nilpotent Lie algebra of maximal class, we have $ \dim L/\gamma_{2}(L)=2 $ and $ \dim \gamma_{i}(L)/\gamma_{i+1}(L)=1 $ for all $ 2 \leq i \leq c $. It is know that the  Frattini subalgebra $ L $ is equal to $ L^{2} $. Hence we may assume that $ L=\langle s, s_{1} \rangle $ such that $ s, s_{1}\notin \gamma_{2}(L) $. On the other hand, by using the proof of $ [5,$ Theorem $1.1] $, we have
\begin{equation*}
\dim(\mathcal{M}(L))=\dim (\mathcal{M}(L/\gamma_{2}(L)))+( \dim(L/\gamma_{2}(L)-1) \dim \gamma_{2}(L) - \sum \limits_{i=2}^{c} \dim \ker(\lambda_{i}).
\end{equation*}
Since $  \dim L/\gamma_{2}(L)=2 $, we have
\begin{equation}
\dim \mathcal{M}(L)= (n-1)- \sum \limits_{i=2}^{n-1} \dim \ker(\lambda_{i}).
\end{equation}
First assume that $ n $ be even. We claim that  $\mathrm{Im}\Psi_{i}\neq 0$  for all odd $ i $ and $2\leq  i \leq n-1$.
Put $ s_{i}=[s_{i-1}, s] $ for all $ i\geq 2 $. We can check that $ s_{i}\in \gamma_{i}(L) \setminus \gamma_{i+1}(L) $ for all $2\leq  i \leq n-1$.  Thun for all $ 2\leq  i \leq n-1 $,
\begin{equation*}
\Psi_{i}(s,s_{1},s,...,s,s_{1},s)=\overline{[s,[s,s_{1},s,...,s,]_{l}] }\otimes \overline{s_{1}}+ \overline{[[s,...,s,s_{1},s]_{r},s]}\otimes \overline{s_{1}}+\overline{w}\otimes \overline{s},
\end{equation*}
for some $w \in L  $. Since $ i $ is odd,
\begin{equation*}
\overline{[s,[s,s_{1},s,...,s,]_{l}] }=\overline{[[s,...,s,s_{1},s]_{r},s]}=\overline{s_{i}}.
\end{equation*}
Hence
$ \Psi_{i}(s,s_{1},s,...,s,s_{1},s)=2\overline{s_{i}}\otimes \overline{s_{1}}+\overline{w}\otimes \overline{s}$. Now, since our field doesn't have  characteristic $2$ and $ s_{i}\in \gamma_{i}(L)\setminus \gamma_{i+1}(L) $,  $2\overline{s_{i}}\otimes \overline{s_{1}} \neq 0$. On the other hand, $ 2\overline{s_{i}}\otimes \overline{s_{1}} $ and $ \overline{w}\otimes \overline{s} $ belong to different direct sum summands $ \gamma_{i}(L)/ \gamma_{i+1}(L) \otimes \langle \overline{s_{1}}\rangle $ and $ \gamma_{i}(L) / \gamma_{i+1}(L) \otimes \langle \overline{s} \rangle $, respectively. Hence  $ 2\overline{s_{i}}\otimes \overline{s_{1}}+\overline{w}\otimes \overline{s} \neq 0$ and so  $\mathrm{Im}\Psi_{i}\neq 0$. By using Corollary $1.2$, we have  $\mathrm{Im}\Psi_{i} \subseteq \ker \lambda_{i}$. Therefore $ \sum \limits_{i=2}^{n-1} \dim \ker(\lambda_{i})\geq\frac{1}{2} (n-2) $. Using  the equality $(2.1)$,
$$ \dim \mathcal{M}(L)=(n-1)-  \sum \limits_{i=2}^{n-1} \dim \ker(\lambda_{i}) \leq (n-1)-(n-2)/2=n/2.    $$
 Let now $ n  $ be odd.  By using Theorem 1.3, we have
$ \dim \mathcal{M}(L)\leq \dim \mathcal{M}(L/Z(L)) + 1$. In this case  $\dim L/Z(L)=n-1$  is even, thus $\dim \mathcal{M}(L/Z(L)) \leq (n-1)/ 2$. Therefore
$\dim \mathcal{M}(L)\leq  (n-1)/ 2+1=(n+1)/2  $ and so $ \dim \mathcal{M}(L) \leq  \lceil(n+1)/2\rceil$. The proof is completed.\\
The following examples show that the upper bound of the Main Theorem  can be obtained by some Lie algebras.
\begin{exa}
 Let $L\cong L(3,4,1,4)=\langle  x_{1},...,x_{4} \mid   [x_{1}, x_{2}]=x_{3}, [x_{1}, x_{3}]=x_{4} \rangle$.  Then $ L $ is a nilpotent Lie algebra of maximal class of dimension $ 4 $. By looking in  $ [1$, Section $4 ] $, we have  $ \dim \mathcal{M}(L(3,4,1,4))=2 $. Hence the bound of the Main Theorem is obtained when $ n $ is even.
\end{exa}
\begin{exa}
Let $L\cong L(7,5,1,7)=\langle x_{1},...,x_{5} \mid [x_{1}, x_{2}]=x_{3}, [x_{1}, x_{3}]=x_{4}, [x_{1}, x_{4}]=x_{5}\rangle $. Then $ L $ is a nilpotent Lie algebra of maximal class of dimension $ 5 $. By using $ [1$, Section $4 ] $, we have  $ \dim \mathcal{M}(L(7,5,1,7))=3 $. Thus the bound of Main Theorem is obtained when $ n $ is odd.
\end{exa}

\end{document}